\documentclass[a4paper,11pt]{article}
\usepackage{graphicx}
\usepackage{subfigure}
\usepackage{epstopdf}
\usepackage{booktabs}
\usepackage{threeparttable}
\usepackage{amssymb, amsmath, amsthm, epsfig}
\numberwithin{equation}{section}
\usepackage{latexsym, enumerate}
\usepackage{eepic}
\usepackage{epic}
\usepackage{color}
\usepackage[colorlinks, linkcolor=blue,anchorcolor=blue,citecolor=blue,urlcolor=blue]{hyperref}

\usepackage{amsmath,textcomp, amssymb, mathrsfs, bm, amsthm, amsfonts}
\usepackage{algpseudocode}
\allowdisplaybreaks[4]
\usepackage{ifpdf}
\usepackage{siunitx}%
\usepackage{chngcntr}

\usepackage{dsfont}
\usepackage{multirow}
\usepackage{bm}
\usepackage{caption}
\allowdisplaybreaks[4]

\theoremstyle{plain}
\newtheorem{lem}{Lemma}[section]
\newtheorem{thm}[lem]{Theorem}

\theoremstyle{definition}
\newtheorem{defn}{Definition}[section]

\theoremstyle{remark}
\newtheorem{rem}{Remark}[section]

\renewcommand{\thefigure}{\thesection.\arabic{figure}}

\newcommand{\ds}{\displaystyle}

\newcommand{\ms}{\medskip}
\begin{document}
\renewcommand{\figurename}{Figure}
\renewcommand{\thesubfigure}{(\alph{subfigure})}
\renewcommand{\thesubtable}{(\alph{subtable})}
\makeatletter
\renewcommand{\p@subfigure}{\thefigure~}

\makeatother
\title{\large\bf Phaseless uniqueness for determining internal source in photo-thermal effect}
\author{
Li-Ping Deng\thanks
{School of Mathematics, Hunan University, Changsha 410082, China.
Email: lipingdeng@hnu.edu.cn}
\and
Hongyu Liu\thanks{Department of Mathematics, City University of Hong Kong, Kowloon, Hong Kong, China.\ \ Email: hongyu.liuip@gmail.com; hongyliu@cityu.edu.hk}
\and
Zhi-Qiang Miao  \thanks
{School of Mathematics, Hunan University, Changsha 410082, China.
Email: zhiqiang\_miao@hnu.edu.cn}
\and
Guang-Hui Zheng\thanks
{School of Mathematics, Hunan University, Changsha 410082, China.
Email: zhenggh2012@hnu.edu.cn; zhgh1980@163.com}
}
\date{}
\maketitle

\begin{center}{\bf ABSTRACT}
\end{center}\smallskip
 The paper investigates an inverse problem of recovering the internal source from external temperature measurements in photo-thermal effect. The photo-thermal effect actually involves two physical processes:  electromagnetic scattering and heat transfer, described by a nonlinear coupled system of Maxwell's equation and the heat transfer equation.
The nonlinear coupling term in the system is represented by the square of the modulus of the electromagnetic (missing the phase information of the electromagnetic field), and the absence of this phase information poses a significant challenge to the reconstruction of the internal source.
In addition, the interaction and mutual influence of multiple physical fields, including electric field, magnetic field and temperature field, add to the complexity involved in the inversion of the internal source.
Based on the potential theory and asymptotic analysis, we prove that the internal source can be uniquely determined up to sign by the external temperature field.
This provides a solid theoretical basis for designing the internal source inversion algorithm and further exploring the theoretical aspects of photo-thermal effect.

\smallskip
{\bf Keywords}: Photo-thermal effect; nonlinear coupling; layer potentials; asymptotic analysis; internal source; uniqueness\\


\section{Introduction}

\subsection{Mathematical setup and statement of main results}
The photo-thermal effect is a fascinating physical phenomenon that occurs when a medium material absorbs light energy, leading to the generation of heat. When light interacts with the surface of a medium material, it is absorbed and subsequently converted into thermal energy, resulting in a rise in the medium material's temperature. This remarkable phenomenon has diverse applications across various fields, such as photothermal energy conversion \cite{Hu2016, Jeong2023, Sasikumar2020}, photothermal therapy \cite{Chu2013, Fernandes2020, Melo2017}, medical imaging \cite{Kim2019, Latterini2014, Sadat2014}, thermoplasmonics \cite{Ammari2018, Baffou2020, Baffou2009, Govorov2007, Gupta2017, Mukherjee2023} and photothermal materials \cite{Fuzil2021, Tee2022, Wu2019}.

The photothermal effect actually involves two fundamental physical processes: electromagnetic scattering and heat conduction. In what follows, let $\Omega \subset \mathbb{R}^{3}$ be a bounded domain occupied by the optical medium material, and such that $\mathbb{R}^{3} \backslash \Omega$ is connected. $Q$ is also a bounded domain that contains the domain $\Omega$. The electric permittivity $\varepsilon_c$ and magnetic permeability $\mu_c$ characterize the optical properties of the medium material $\Omega$. The thermal properties of $\Omega$ are described by mass density $\rho_c$, specific heat $\alpha_c$, and thermal conductivity $\gamma_c$. On the other hand, the optical and thermal properties of homogeneous medium material in $\mathbb{R}^{3} \backslash \Omega$ are characterized by electric permittivity $\varepsilon_m$, magnetic permeability $\mu_m$, mass density $\rho_m$, specific heat $\alpha_m$, and thermal conductivity $\gamma_m$, respectively. Moreover, we assume that the medium material is nonmagnetic, i.e., $\mu_c=\mu_m$, and the permittivity of the medium material is given by the Drude Model:
\begin{equation}\label{Drude}
\varepsilon_c(\omega)=1-\frac{\omega_p^2}{\omega(\omega+i\tau)},
\end{equation}
where $i:=\sqrt{-1}$, $\omega_p$ is the electric plasma frequency, and $\tau$  is the damping coefficient.

Set
\begin{align*}
\varepsilon=\varepsilon_c \chi(\Omega)+\varepsilon_m \chi(\mathbb{R}^{3} \backslash \overline {\Omega}),
        \quad \rho=\rho_c \chi(\Omega)+\rho_m \chi(\mathbb{R}^{3} \backslash \overline {\Omega}),\\
\alpha=\alpha_c \chi(\Omega)+\alpha_m \chi(\mathbb{R}^{3} \backslash \overline {\Omega}) ,
        \quad \gamma=\gamma_c \chi(\Omega)+\gamma_m \chi(\mathbb{R}^{3} \backslash \overline {\Omega}) ,
\end{align*}
where $\chi(\Omega)$ denotes the characteristic function of $\Omega$.
Let $E(x,\omega)$, $H(x,\omega)$, $(x,\omega)\in \mathbb{R}^3\times\mathbb{R}_+$ denote the time-harmonic electric field, magnetic field, and $v(x)$, $x\in \mathbb{R}^3$ denotes steady-state temperature field, respectively. The photo-thermal effect is modeled by the following electromagnetic scattering and heat transfer coupling system:
\begin{align}\label{model2}
\begin{cases}
\ds   \nabla \times E-i\omega\mu H=0  &\quad \text{in} \quad \mathbb{R}^{3},\\
\ds   \nabla \times H+i \omega\varepsilon E=J  &\quad \text{in} \quad \mathbb{R}^{3},\\
\ds   \nabla \cdot H=0,\quad \nabla \cdot E=\frac{\rho}{\varepsilon}  &\quad \text{in} \quad \mathbb{R}^{3},\\
\ds   E=E^i+E^s, \quad  H=H^i+H^s  &\quad \text{in} \quad \mathbb{R}^{3},\\
\ds   \lim\limits_{|x|\rightarrow \infty}|x|(\sqrt{\mu}H^sx\hat{x}-\sqrt{\varepsilon}E^s)=0  &\quad \hat{x}:=x/|x|,\\
\ds   \nabla \cdot(\gamma\nabla v)=\frac{\omega}{8\pi}  Im(\varepsilon)|E|^2  &\quad \text{in} \quad \mathbb{R}^{3},\\
\ds   v-V(x)=O\Big(\frac{1}{|x|}\Big)  &\quad |x|\rightarrow\infty,\\
\end{cases}
\end{align}
where $V(x)$ is the background temperature field, $J(x,\omega)$ represents the electric current. The electromagnetic scattering process is triggered by a specified time-harmonic incidence field $(E^i, H^i)$, which satisfy
\begin{equation*}
\nabla \times E^i-i\omega H^i=0,\quad \nabla \times H^i+i \omega\varepsilon_m E^i=0 \qquad \text{in}\quad \mathbb{R}^{3},\\
\end{equation*}
and the first limit in \eqref{model2} is known as the Silver-M$\ddot{\rm u}$ller radiation condition, which is employed to ensure the well-posedness of the scattering problem. We would like to emphasize that the electromagnetic scattering process and the heat transfer process are coupled through the nonlinear term $\frac{\omega}{8\pi}  Im(\varepsilon)|E|^2 $, the electric field with phaseless information acts as the heat source to generate the temperature field. For the electromagnetic scattering and heat transfer coupling system, we are particularly focused on the following inverse problem:
\begin{equation}\label{inverse problem1}
v|_{(x,\omega)}\in \partial Q\times R_+\longrightarrow J(x,\omega).\\
\end{equation}
That is to say, the objective is to recover the internal source by measuring the external temperature field. The following analysis shows that the absence of phase information in the coupling term mentioned above causes considerable difficulties for the inversion of the internal source.
\par Next, we will study a particular scenario in which the medium material $\Omega$ is characterized by an infinitely long cylinder with the axis coincident with the $x_3$-axis for $x=(x_j)^3_{j=1}\in \Omega$,
 and we take the form of an infinitely cylinder with the transverse profile $D$ along the $x_3$-axis.
 Additionally, we define $B_R:=B_R(0)$ as the disk centered at the origin with radius $R \in \mathbb R_+$, which contains the transverse profile $D$.
\par Transverse magnetic polarization is a type of electromagnetic wave polarization commonly used in the study of wave propagation and scattering, which is often employed in various applications, including optics, microwave engineering, and electromagnetic modeling, due to its distinct characteristics and convenient mathematical formulations.
Under the transverse-magnetic (TM) polarization, which refers to \cite{Bondarenko2016,Cao2020}, that is,
\begin{align*}
E^i=
\left[
\begin{array}{c}
  0  \\
  0  \\
  u^i(x_1,x_2)
\end{array}
\right],\quad
H^i=
\left[
\begin{array}{c}
  H_1^i(x_1,x_2)  \\
  H_2^i(x_1,x_2)  \\
  0
\end{array}
\right],\\
\end{align*}
and
\begin{align*}
E=
\left[
\begin{array}{c}
  0  \\
  0  \\
  u(x_1,x_2)
\end{array}
\right],\quad
J=
\left[
\begin{array}{c}
  0  \\
  0  \\
  f(x_1,x_2)
\end{array}
\right],\quad
H=
\left[
\begin{array}{c}
  H_1(x_1,x_2)  \\
  H_2(x_1,x_2)  \\
  0
\end{array}
\right],\\
\end{align*}
therefore, through straightforward calculations, the time-harmonic model \eqref{model2} can be reduced to
\begin{align}\label{problem model}
\begin{cases}
   \Delta u +\omega^{2} \varepsilon_c(\omega)u =-i\omega f(x)  &\quad\text{in} \ \ D , \\
   \Delta u +\omega^{2} \varepsilon_m u =0  &\quad\text{in} \ \ \mathbb{R}^{2} \backslash \overline {D}, \\
   u_{+}=u_{-},\quad \frac{\partial u}{\partial\nu}\big{|}_{+}=\frac{\partial u}{\partial\nu}\big{|}_{-} &\quad \text{on}\ \ \partial D,\\
    \lim\limits_{r\rightarrow\infty}\sqrt{r}(\partial_r u-iku)=0 &\quad r=|x|,\\
   \nabla(\gamma_c \cdot \nabla v)=\frac{\omega}{2\pi} Im\left(\varepsilon_{c}(\omega)\right) |u|^2 & \quad \text{in} \ D, \\
   \nabla(\gamma_m \cdot \nabla v)=0  &\quad \text{in} \ \ \mathbb{R}^{2} \backslash \overline {D}, \\
   v_+=v_-,\quad  \gamma_m \frac{\partial v}{\partial \nu}\big{|}_{+}=\gamma_c \frac{\partial v}{\partial \nu}\big{|}_{-}  &\quad \text{on}\ \ \partial D,\\
   v-V(x)=O\big(\frac{1}{|x|}\big) &\quad \lvert x \rvert \rightarrow \infty,
\end{cases}
\end{align}
where $\frac{\partial}{\partial \nu}$ signifies the outward normal derivative, and we use the notation $\frac{\partial}{\partial \nu}|_{\pm}$ indicating
$$ \frac{\partial v}{\partial \nu}\bigg|_{\pm}(x)=\lim_{t\rightarrow 0^+}\nabla v(x\pm t\nu(x))\cdot \nu(x)$$
with $\nu$ as the outward unit normal vector to $\partial D$.
The variable $f(x)$  denotes the internal source and  $f(x)\in L^\infty(\mathbb R^2)$ with supp$(f)\in D \subset B_R$. Furthermore, we perform dimensionless normalization on parameters $\varepsilon_m, \mu_m , \gamma_m$, that is $\varepsilon_m=\mu_m =\gamma_m=1$. It is emphasized that we assume the medium material is nonmagnetic, that is, $\mu_c=\mu_m=1$.
\par After considering the transverse magnetic (TM) polarization in time-harmonic and steady-state system, the inverse problem \eqref{inverse problem1} that we are interested in is transformed into the following one:
\begin{equation}\label{inverse problem2}
v|_{(x,\omega)}\in \partial B_R\times R_+\longrightarrow f(x).\\
\end{equation}
\par For simplicity, we further introduce the following boundary measurement operator:
\begin{equation}\label{boundary measurement}
\Lambda_{f}(x,\omega)=v(x,\omega),\qquad (x,\omega)\in \partial B_R\times \mathbb{R_+}.\\
\end{equation}
Then, the inverse problem \eqref{inverse problem2} that we are concerned with is to determine the unknown internal source $f(x)$ by knowledge of $\Lambda_{f}(x,\omega),(x,\omega)\in \partial B_R\times \mathbb{R_+}$, that is
\begin{equation}\label{inverse problem3}
\Lambda_{f}|_{(x,\omega)}\in \partial B_R\times R_+\longrightarrow f(x).\\
\end{equation}

To state our main results for determining the internal source in photo-thermal effect, we first give the
following definitions that describe the a priori properties of the internal source.

\begin{defn}\label{adm}
If internal sources $f_1$, $f_2$ satisfy $f_1\pm f_2=h(x)$, and $h(x)$ is harmonic function in $\mathbb R^2$, we say that the two internal sources $f_1$ and $f_2$ differ by a harmonic part. Furthermore, if $f_1\pm f_2=l(x)$, and $l(x)$ is biharmonic function in $\mathbb R^2$, we say that the two internal sources differ by a biharmonic part. Clearly, two internal sources that differ by a harmonic part also differ by a biharmonic part. Conversely, the opposite statement does not hold.
\end{defn}

\begin{defn}\label{adm}
We call the two internal sources $f_1$ and $f_2$ invariant along a direction vector $\eta$ in $\mathbb R^2$, if $\nabla(f_1\pm f_2)\cdot\eta=0$.
\end{defn}

\par Now, we present the main result of the inverse problem \eqref{inverse problem3}. Indeed, we establish the phaseless uniqueness result for the inverse problem \eqref{inverse problem3} as follows.

\begin{thm}\label{Th}
Consider the coupling system \eqref{problem model} association with two internal sources $f_1$ and $f_2$. Assume that the internal sources $f_1$ and $f_2$  are supported in $D$, and satisfy either of the following conditions:\\
(i)\  $f_1$ and $f_2$ differ by a harmonic part or biharmonic part;\\
(ii)\ $f_1$ and $f_2$ are invariant along a direction vector $\eta$ in $\mathbb R^2$.\\
If
\begin{equation}\label{equation1}
\Lambda_{f_1}(x,\omega)=\Lambda_{f_2}(x,\omega),\qquad (x,\omega)\in \partial B_R\times \mathbb{R_+},\\
\end{equation}
we have
\begin{equation*}
|f_1(x)|=| f_2(x)|,\qquad x\in D.\\
\end{equation*}
\end{thm}

\begin{rem}
From Theorem \ref{Th}, the uniqueness holds (up to sign), and the essential reason is the phaseless coupling between the electromagnetic field and the temperature field. Obviously, if the internal source remains positive or negative, the source is uniquely determined.
\end{rem}

\subsection{Related results and discussion}

The issue of inverse source scattering, specifically the uniqueness and non-uniqueness of solutions, is of great interest and significance in various fields of science and engineering. The uniqueness of inverse source scattering problems is influenced by several factors, including the nature of the scatterer, the type of incident waves, and the available measurement configurations. However, non-uniqueness poses a challenge to the reconstruction process as it introduces uncertainties and complicates the interpretation of the measured data. In general, there are a number of papers claiming that the inverse-source and scattering problems do not have unique solutions \cite{Bleistein1977,Devaney1982,Hauer2005}, additional constraints need to be imposed to obtain a unique solution.

For the acoustic waves or the Helmholtz equation, there are a lot of mathematical results available.
The reference \cite{Bao2010} is concerned with an inverse source problem, determining the source from measurements of the radiated fields away at multiple frequencies. In \cite{Liu2015}, the authors established the uniqueness of identifying both the acoustic density and the internal source for the Helmholtz equation (or time domain wave equation) through passive measurements, a scenario that arises in thermoacoustic and photoacoustic tomography. In \cite{ji2019}, the authors propose a stable phase retrieval method and prove several results regarding uniqueness using phaseless far-field data for inverse source scattering problems.
The problem is to determine both the source and its velocity in a medium based on measurements of the solution of the wave equation on the boundary in \cite{Knox2020}.
In \cite{Zhangd2023(1)}, the authors demonstrate that it is feasible  to simultaneously recover the positions of the sources and the obstacle using the scattered Cauchy data of the time-harmonic acoustic field.

In the case of electromagnetic waves,
the reference \cite{Deng2019} considers the unique recovery of an unknown electric current and surrounding parameters of a medium from boundary measurements, which is particularly relevant in the field of brain imaging.
In \cite{ammari2020}, the authors offer a rigorous mathematical analysis of the anticipated superresolution phenomenon in time-reversal imaging of electromagnetic (EM) radiating sources situated within a high-contrast medium. In \cite{ji201901}, the authors investigate the uniqueness and numerical algorithm (such as phase retrieval, shape reconstruction methods) for inverse electromagnetic source scattering problems with multifrequency sparse phased or phaseless far field data.

In recent years, there has been an increasing interest in the study of elastic wave scattering problems.
In \cite{Bao2020}, the authors focus on the inverse source problems related to the time-harmonic elastic wave equations. Their study indicates that the external force and electric current density can be identified based on the boundary measurement of the radiated wave field.
In \cite{Li2020}, the authors provide a demonstration that the mean amplitude of the scattering field over the specified frequency range obtained from a single instance of the random source is sufficient to determine the principle symbol of the covariance operator.
In \cite{Chang2023}, the objective of this study is to reconstruct a rigid obstacle and the excitation sources simultaneously using near-field measurements.

Moreover, there is a lot of  literature available on relevant stability estimates and algorithmic research\cite{Bao2020,Bao2010,Chang2021,Chang2023,Lipeijun2020,Zhangd2023(1),Zhangd2023(2)}.



In this work, we utilize an electromagnetic scattering and heat transfer coupling system to simulate the photo-thermal effect.
Furthermore, we study the phaseless uniqueness of internal source determination in this coupling physical process, which makes the topic very challenging. First, it involves the nonlinear coupling of two partial differential equations. Multiple physical processes and unknown parameters are interdependent and interact with one another, which complicates our analysis significantly. Secondly, the coupling term is the phaseless, and the absence of phase information significantly exacerbates the ill-posedness of the inverse internal source problem. To overcome these difficulties and challenges, we use the layer potential theory, asymptotic analysis techniques to solve the inverse problem. For the first part of the coupled system, i.e., the electromagnetic scattering process, we mainly utilize low frequency asymptotic expansion to analyse the electromagnetic field . The second part is the instantaneous heat process, which is described by the steady-state heat equation and we solve it using the layer potential theory. Finally, we derive the integral identities involving unknown parameter, which allow us to establish the phaseless uniqueness of recovering the internal source.

\par The rest of the paper is organized as follows. In Section \ref{sec:pk}, we provide the preliminary knowledge. In Section \ref{sec:Asymptotic analysis}, we utilize asymptotic analysis techniques and layer potential theory to dissect electromagnetic scattering. The steady-state heat conduction process is analysed in Section \ref{sec:Asymptotic analysis-heat}. The proof of the phaseless uniqueness for the coupled system \eqref{problem model} is presented in Section \ref{sec:proof}.

\section{Preliminary knowledge}\label{sec:pk}
\noindent
The purpose of this chapter is to provide readers with the essential background knowledge necessary to comprehend the research topic addressed in the paper. The investigation undertaken in this study mainly relies on the asymptotic analysis techniques and layer potential theory.
\par Let $\Phi(\cdot,k)$ be the Green function for the Helmholtz operator $\Delta+k^2$ satisfying the Sommerfeld radiation. In two dimension, $\Phi(x,k)$ is given by
$$\Phi(x,k):=-\frac{i}{4}H_0^{(1)}(k|x|),$$
where $H_0^{(1)}$ is the Hankel function of first and order 0. We further let $\Phi(x,y,k):=\Phi(x-y,k)$. For $k\ll 1$ we have the asymptotic expansion
\begin{equation}\label{expansion1}
\Phi(x-y, k)
=\Phi_0(x,y)+N+\frac{1}{2\pi}\ln k -\frac{i}{4}+\sum_{j=1}^\infty (b_j \ln k \lvert x-y \rvert +c_j)(k \lvert x-y \rvert)^{2j}.\\
\end{equation}
where
\begin{align*}
\begin{cases}
\ds \Phi_0(x,y)=\frac{1}{2\pi} \ln{\lvert x-y \rvert},\\
\ds N=\frac{1}{2\pi}(\gamma_e-\ln2),\\
\ds \gamma_e=\lim\limits_{p\rightarrow \infty} \big(\sum\limits_{m=1}^p \frac{1}{m}-\ln p \big) ,\\
\ds b_j=\frac{(-1)^j}{2\pi} \frac{1}{2^{2j}(j!)^2},\\
\ds c_j=-b_j(\gamma_e-\ln 2-\frac{i\pi}{2}-\sum\limits_{n=1}^j\frac{1}{n}),\\
\end{cases}
\end{align*}
where $\gamma_e$ is the Euler constant and $N$ is a constant(see \cite{Ammari2018}). It is known that $\Phi_0(x,y)$  is the fundamental solution of the Laplace operator $\Delta$.

\par In the sequel, we define the following single and double layer
potentials by $\Phi_0$, for a bounded domain $D$ is of class $C^{1,\alpha}$ for some $0<\alpha<1$ and the density function $\varphi(x)\in L^2(\partial D)$ \cite{Colton2019}:
\begin{equation*}
\mathcal S_D: H^{-\frac{1}{2}}(\partial D)\rightarrow  H^{\frac{1}{2}}(\partial D),\quad
 \mathcal S_D[\varphi](x):=\displaystyle\int_{\partial D} \Phi_0(x,y)\varphi(y)d\sigma(y),\\
\end{equation*}

Next, we will introduce a lemma known as the jump relation,  which is a crucial auxiliary tool in the analysis of potential theory \cite{Ammari2010}.

\begin{lem}[jump relation]\label{lem:jump relation}
For the density function ~$\varphi(x)\in L^2(\partial D)$~, it is well known that
\begin{equation}
\label{eq:jump relation}
  \frac{\partial \mathcal S_D[\varphi]}{\partial \nu}\big|_{\pm}(x)=\left(\pm \frac{1}{2}I+ \mathcal K_D^* \right)[\varphi](x),\quad a.e. \quad x\in \partial D.\\
\end{equation}
where
the operator $\mathcal K_D^*$ is the $L^2$-adjoint of $\mathcal K_D$ defined by
\begin{equation*}
  \mathcal K_D^*: H^{-\frac{1}{2}}(\partial D)\rightarrow  H^{-\frac{1}{2}}(\partial D),\quad
  \mathcal K_D^*[\varphi](x):=p.v.\displaystyle\int_{\partial D} \frac{\partial\Phi_0(x,y)}{\partial\nu(x)}\varphi(y)d\sigma(y),
\end{equation*}
here p.v. stands for the Cauchy principal value.
\end{lem}

Furthermore, we also need to provide the definition of volume potential, which is a type of potential field represented by a density function distributed within a region $D$. For a function $\varphi(x)\in L^\infty(\mathbb R^2)$ supported in $D$, the volume potential $\mathcal N_D$ can be defined as follows \cite{Colton2019}:
\begin{equation*}
 \mathcal N_D: L^2(D)\rightarrow L^2(D),\quad
 \mathcal N_D[\varphi](x):=\int_D\Phi_0(x,y)\varphi(y)dy.
\end{equation*}

\par In this article, the permittivity of medium material is given by the Drude Model \eqref{Drude}, and it can be rewritten as,
\begin{equation*}
\label{Drude2}
\varepsilon_c(\omega)=1-\frac{\omega_p^2}{\omega(\omega+i\tau)}=1-\frac{\omega_p^2}{\omega^2+\tau^2}
+i\frac{\omega_p^2 \tau}{\omega(\omega^2+\tau^2)}.
\end{equation*}
Then we have the asymptotic expansions of the real and imaginary parts of $\varepsilon_c(\omega)$ as $\omega\rightarrow 0^+$:
\begin{equation*}
Re\big(\varepsilon_c(\omega)\big)=1-\frac{\omega_p^2}{\tau^2}+\omega^2\frac{\omega_p^2}{\tau^4}+O(\omega^4),
\end{equation*}
and
\begin{equation*}
Im\big(\varepsilon_c(\omega)\big)=\frac{1}{\omega}\frac{\omega_p^2}{\tau}-\omega\frac{\omega_p^2}{\tau^3}+O(\omega^3).
\end{equation*}
Consequently, $\varepsilon_c(\omega)$ can be expressed by
\begin{equation}\label{expansion Drude}
\varepsilon_c(\omega)=\frac{1}{\omega}\frac{i\omega_p^2}{\tau}+\left(1-\frac{\omega_p^2}{\tau^2}\right)-\omega\frac{i\omega_p^2}{\tau^3}
+\omega^2\frac{\omega_p^2}{\tau^4}+O(\omega^3).
\end{equation}
Through the aforementioned asymptotic expansion analysis, we can approximate the real and imaginary parts of the Drude model at a low frequency.

\par In the rest of this section, we will introduce a series of definitions along with a critical lemma, which provide a theoretical foundation for subsequent analysis \cite{Ammari2018}.
\begin{defn}\label{def4.1}
Let $\mathcal C=\left\{\varphi \in H^{-\frac{1}{2}}(\partial D); \exists\  \alpha \in\mathbb C, S_D[\varphi]=\alpha\right\}.$
We call $\varphi_0$ the unique element of $\mathcal C$ such that $\displaystyle\int_{\partial D}\varphi_0 d\sigma=1. $
\end{defn}

\begin{defn}\label{def:S}
Let the linear operator $\tilde{S}_D$ satisfy
\begin{equation*}
\tilde{\mathcal S}_D:\  H^{-\frac{1}{2}}(\partial D)\rightarrow  H^{\frac{1}{2}}(\partial D),\\
\ \varphi\rightarrow
\begin{cases}
\mathcal S_D[\varphi] &\qquad if\ (\varphi,1)_{-\frac{1}{2},\frac{1}{2}}=0,\\
-1 &\qquad if\  \varphi_0=\varphi.
\end{cases}
\end{equation*}
where $(\cdot\ ,\ \cdot)$ is the dual product with the space $H^{-\frac{1}{2}}(\partial D)$ and $H^{\frac{1}{2}}(\partial D)$.
\end{defn}

\begin{rem}
In fact, single layer potential operator $S_D$ is not invertible in $\mathbb R^2$, therefore we define the invertible operator $\tilde{S}_D$ to substitute for $S_D$. For the sake of symbol unification, we will continue to denote $S_D$ as $\tilde{S}_D$.
\end{rem}

\begin{lem}\label{lem4.1}
Let $g\in H^1(D)$ be such that $\Delta g=f$ with $f\in L^2(D)$. Then, in $H^*(\partial D)$,
\begin{equation*}
\left(\frac{1}{2}I-K_D^*\right)\mathcal S_D^{-1}[g]=-\frac{\partial g}{\partial \nu}+T_f,\\
\end{equation*}
where $T_f\in H^*(\partial D)$ and $\|T_f\|_{H^*(\partial D)}\leq C\|f\|_{L^2(D)}$ for a constant C.\\
\par Moreover, if $g\in H^1_{loc}(\mathbb{R}^{2}),\Delta g=0$ in $\mathbb{R}^{2} \backslash \overline {D}$ and $\lim\limits_{|x|\rightarrow \infty} g(x)=0$, then
\begin{equation*}
T_f=c_f\varphi_0+\mathcal S_D^{-1}[g]\\
\end{equation*}
with
$$ c_f=\int_Df(x)dx-\int_{\partial D}\mathcal S_D^{-1}[g](y)d\sigma(y), $$
where $\varphi_0$ is given in Definition \ref{def4.1}. Here, by an abuse of notation, we still denote the trace of g on $\partial D$.\\
\end{lem}

\par In addition, we also introduce a lemma to supplement our analysis. This lemma serves to provide additional support and clarification for the arguments presented in the paper\cite{Ammari2010,Ammari2004}.
\begin{lem}\label{lem:inverse}
If $S_D$ and $(\lambda I-K_D^*)$ is invertible, then we have the inverse
\begin{equation*}
\left(
\begin{array}{cc}
  \mathcal S_D & -\mathcal S_D  \ms\\
  \frac{1}{2}I+ \mathcal K_D^* & -\gamma_c\left(-\frac{1}{2}I+\mathcal K_D^*\right)
\end{array}
\right)^{-1}
=\frac{1}{\gamma_c-1}
\left(
\begin{array}{cc}
  \gamma_c\left(\lambda I-\mathcal K_D^*\right)^{-1}\left(\frac{1}{2}I-\mathcal K_D^*\right)\mathcal S_D^{-1} & \left(\lambda I-\mathcal K_D^*\right)^{-1} \ms  \\
  -\left(\lambda I-\mathcal K_D^*\right)^{-1}\left(\frac{1}{2}I+\mathcal K_D^*\right)\mathcal S_D^{-1} & \left(\lambda I-\mathcal K_D^*\right)^{-1}
\end{array}
\right)
\end{equation*}
where
\begin{equation}\label{lambda}
\lambda=\frac{\gamma_c+1}{2(\gamma_c-1)}.
\end{equation}
\end{lem}
It is noted that $\lambda$ in \eqref{lambda} is not a Dirichlet eigenvalue for $-\Delta$ on $D$. Hence $(\lambda I-\mathcal K_D^*):L^2(D)\rightarrow L^2(D)$ is injective.

\section{Asymptotic analysis for electromagnetic scattering }\label{sec:Asymptotic analysis}
\par In the first process, we mainly use the Lippmann-Schwinger integral equation and the asymptotic expansion to analyze the electric field.

\par Due to the inherent losses of the medium material, we further assume that
\begin{equation*}
Im\left(\varepsilon_c(\omega) \right)\geq 0.\\
\end{equation*}

\begin{lem}
Assume $f(x)\in L^\infty(\mathbb R^2)$, and supp$(f)\in D \subset B_R$, the nonlinear coupled system \eqref{problem model} is uniquely solvable.
\end{lem}
\begin{proof}
If the internal source $f(x)$ is known, from Lippmann-Schwinger integral equation \cite{Colton2019} and classical elliptic equation theory \cite{McLean2000}, the conclusion is obvious.
\end{proof}

\begin{lem}\label{lem3.2}
Let $u(x,\omega)\in H_{loc}^1(\mathbb R^2)$ be the solution to \eqref{problem model}.
Then $u(x,\omega)$ is uniquely given by the following integral equation
\begin{equation}\label{LS}
u(x,\omega)=\omega^2\int_{\mathbb R^2}(1- \varepsilon) u(y) \Phi(x,y)dy-i\omega
\int_{\mathbb R^2}\Phi(x,y)f(y)dy,\quad x\in\mathbb{R}^2. \\
\end{equation}
Furthermore, as $\omega\rightarrow +0$, we have
\begin{equation}\label{equation u}
\begin{split}
u(x,\omega)=
&-\omega\ln\omega \frac{i}{2\pi}\mathcal T_{1}[f]-\omega\left[\left(\frac{1}{4}+iN\right)\mathcal T_{1}[f]+i\mathcal N_{D}[f](x)\right]
-\omega^2\ln^2\omega\frac{\omega_p^2}{4\pi^2\tau}|D|\cdot\mathcal T_{1}[f]+O(\omega^2\ln\omega), \\
\end{split}
\end{equation}
where $|D|$ denotes the area of domain $D$, and $\mathcal T_{1}[f]:=\displaystyle\int_Df(y)dy$.
\end{lem}

\begin{proof}
Since the integral equation \eqref{LS} is of Lippmann-Schwinger type,
it is uniquely solvable for $u\in C(\mathbb{R}^2)$(see \cite{Colton2019}).
For $x \in \mathbb{R}^2 $, we have
\begin{equation*}
\begin{split}
u(x,\omega)=&\omega^2\int_D [1- \varepsilon_c(\omega)] u(y) \Phi(x,y)dy-i\omega\int_D \Phi(x,y) f(y)dy, \\
\end{split}
\end{equation*}
By using the asymptotic expansions \eqref{expansion1} and \eqref{expansion Drude},  and  keeping the highest term to $O(\omega\ln\omega)$,  together with
straightforward calculations,  we can see that
\begin{align*}
u(x,\omega)
=&-\omega\ln\omega\frac{i\omega_p^2}{2\pi \tau}\int_Du(y)dy-\omega\frac{i\omega_p^2}{\tau}\int_Du(y)\left[\Phi_0(x,y)+N-\frac{i}{4} \right]dy-\omega\ln\omega\frac{i}{2\pi}\int_Df(y)dy\\
&-i\omega\int_D\Phi_0(x,y)f(y)dy-i\omega\left( N-\frac{i}{4} \right)\int_Df(y)dy+O(\omega^2\ln\omega)\\
=&-\omega\ln\omega \frac{i}{2\pi}\mathcal T_{1}[f]-\omega \left(\frac{1}{4}+iN\right)\mathcal T_{1}[f]-i\omega \mathcal N_{D}[f](x)-\mathcal T_{\omega}[u](x)+O(\omega^2\ln\omega),
\end{align*}
where \begin{equation}\label{def formula1}
\mathcal T_{\omega}[u](x):=\omega\ln\omega\frac{i\omega_p^2}{2\pi\tau}\mathcal T_{1}[u]+\omega\left[\frac{i\omega_p^2}{\tau}\mathcal N_{D}[u](x)
+\frac{i\omega_p^2}{\tau}\left(N-\frac{i}{4}\right)\mathcal T_{1}[u]\right].
\end{equation}
Obviously, the above solution is an implicit expression of $u (x)$, and we have
\begin{align*}
(I+\mathcal T_\omega)[u](x,\omega)
=-\omega\ln\omega \frac{i}{2\pi}\mathcal T_{1}[f]-\omega \left(\frac{1}{4}+iN\right)\mathcal T_{1}[f]-i\omega \mathcal N_{D}[f](x)+O(\omega^2\ln\omega).
\end{align*}
To see this, notice first that $\omega\ln\omega\rightarrow 0$ as $\omega\rightarrow 0$, so $\mathcal T_\omega[u](x)\rightarrow 0$ as $\omega\rightarrow 0$, then we deduce that
\begin{equation*}
\begin{split}
(I+\mathcal T_\omega)^{-1}=&I-\mathcal T_\omega+\mathcal T_\omega^2-\mathcal T_\omega^3+\cdots \\
=&I-\omega\ln\omega\frac{i\omega_p^2}{2\pi\tau}\mathcal T_1[\cdot]
-\omega\left[\frac{i\omega_p^2}{\tau}\mathcal N_{D}[\cdot](x)
+\frac{i\omega_p^2}{\tau}\left(N-\frac{i}{4}\right)\mathcal T_1[\cdot]\right]+O(\omega^2\ln^2\omega).
\end{split}
\end{equation*}
Plugging this identity into the following equation
\begin{align*}
u(x,\omega)=&(I+\mathcal T_\omega)^{-1}\left[-\omega\ln\omega \frac{i}{2\pi}\mathcal T_{1}[f]-\omega \left(\frac{1}{4}+iN\right)\mathcal T_{1}[f]-i\omega \mathcal N_{D}[f](x)+O(\omega^2\ln\omega)\right].\\
\end{align*}
Consequently, as $\omega \rightarrow +0$,
\begin{align*}
u(x,\omega)=& \left\{I-\omega\ln\omega\frac{i\omega_p^2}{2\pi\tau}\mathcal T_1[\cdot] - \omega\left[\frac{i\omega_p^2}{\tau}\mathcal N_{D}[\cdot](x)
+\frac{i\omega_p^2}{\tau}\left(N-\frac{i}{4}\right)\mathcal T_1[\cdot]\right]+O(\omega^2\ln^2\omega) \right\} \\
& \cdot\left[-\omega\ln\omega \frac{i}{2\pi}\mathcal T_{1}[f]-\omega \left(\frac{1}{4}+iN\right)\mathcal T_{1}[f]-i\omega \mathcal N_{D}[f](x)+O(\omega^2\ln\omega)\right]\\
=&-\omega\ln\omega \frac{i}{2\pi}\mathcal T_{1}[f]-\omega\left(\frac{1}{4}+iN\right)\mathcal T_{1}[f]-i\omega \mathcal N_{D}[f](x)-\omega^2\ln^2\omega\frac{\omega_p^2}{4\pi^2\tau}\mathcal T_1[\mathcal T_{1}[f]]+O(\omega^2\ln\omega)\\
=&-\omega\ln\omega \frac{i}{2\pi}\mathcal T_{1}[f]-\omega\left[\left(\frac{1}{4}+iN\right)\mathcal T_{1}[f]+i\mathcal N_{D}[f](x)\right]
-\omega^2\ln^2\omega\frac{\omega_p^2}{4\pi^2\tau}|D|\cdot\mathcal T_{1}[f]+O(\omega^2\ln\omega),
\end{align*}
then we derive the equation \eqref{equation u}.
\end{proof}

According to the equation \eqref{equation u}, we easily have the real part and the imaginary part of $u$,  as $\omega \rightarrow +0$
\begin{equation*}
\begin{split}
Re(u)=&-\frac{1}{4}\omega \mathcal T_{1}[f]-\omega^2\ln^2\omega\frac{\omega_p^2}{4\pi^2\tau}|D|\cdot\mathcal T_{1}[f]+O(\omega^2\ln\omega),
\end{split}
\end{equation*}
and
\begin{equation*}
\begin{split}
Im(u)=&-\left\{\omega\ln\omega\frac{1}{2\pi}\mathcal T_{1}[f]+\omega\left[N\mathcal T_{1}[f]+\mathcal N_{D}[f](x)\right]\right\}+O(\omega^2\ln\omega).\\
\end{split}
\end{equation*}

\section{Asymptotic analysis for heat transfer}\label{sec:Asymptotic analysis-heat}
Since  $\gamma_m = 1$, the steady-state heat transfer problem can be simplified to
\begin{align}\label{problem2.1}
\begin{cases}
\ds     \Delta v=\frac{\omega}{2\pi\gamma_c} Im(\varepsilon_{c}) |u|^2  &\quad\text{in} \ \ D , \\
\ds     \Delta v=0 &\quad \text{in} \ \ \mathbb{R}^{2} \backslash \overline {D}, \\
\ds     v_+-v_-=0  &\quad \text{on} \ \ \partial D,\ms\\
\ds     \frac{\partial v}{\partial \nu}\big{|}_{+}-\gamma_c \frac{\partial v}{\partial \nu}\big{|}_{-}=0 &\quad \text{on} \ \ \partial D,\\
\ds     v-V(x)=O\Big(\frac{1}{\lvert x \rvert}\Big) &\quad \lvert x \rvert \rightarrow \infty.
\end{cases}
\end{align}

\par Before  solving the temperature field $v (x, \omega)$, we need to calculate the right term of \eqref{problem2.1}. We first compute the square of the electric field $u$ as
\begin{align*}
\left|u(x,\omega)\right|^2=&Re^2(u)+Im^2(u) \\
=&\left[-\frac{1}{4}\omega \mathcal T_{1}[f]-\omega^2\ln^2\omega\frac{\omega_p^2}{4\pi^2\tau}|D|\cdot\mathcal T_{1}[f] \right]^2 +\left[\omega\ln\omega\frac{1}{2\pi}\mathcal T_{1}[f]+\omega\left(N\mathcal T_{1}[f]+\mathcal N_{D}[f](x)\right) \right]^2 +O(\omega^3\ln\omega)\\
=&\omega^2\ln^2\omega\frac{1}{4\pi^2}(\mathcal T_{1}[f])^2+\omega^2\ln\omega\frac{1}{\pi}\big(\mathcal T_{1}[f]\mathcal N_{D}[f](x)+N(\mathcal T_{1}[f])^2\big)+\omega^2(\mathcal N_{D}[f])^2(x)\\
&+\omega^2\left(N^2+\frac{1}{16}\right)(\mathcal T_{1}[f])^2+2\omega^2N\mathcal T_{1}[f]\mathcal N_{D}[f](x)+O(\omega^3\ln\omega).
\end{align*}
Hence
\begin{align*}
\omega \frac{Im(\varepsilon_c)}{2\pi\gamma_c}\left|u\right|^2
=&\omega^2\ln^2\omega\frac{\omega_p^2}{8\pi^3\gamma_c\tau}(\mathcal T_{1}[f])^2+\omega^2\ln\omega\left\{\frac{\omega_p^2}{2\pi^2\gamma_c\tau}\left[\mathcal T_{1}[f]\mathcal N_{D}[f](x)+N(\mathcal T_{1}[f])^2\right]\right\}\\
&+\omega^2\frac{\omega_p^2}{2\pi\gamma_c\tau}\left[ (\mathcal N_{D}[f])^2(x)+\left(N^2+\frac{1}{16}\right)(\mathcal T_{1}[f])^2+2N\mathcal T_{1}[f]\mathcal N_{D}[f](x)\right]+O(\omega^3\ln\omega).
\end{align*}

To simplify subsequent calculations, it is convenient to set
\begin{equation*}
\begin{split}
Q(x,\omega)
=&\omega^2\ln^2\omega Q_1(x)+\omega^2\ln\omega Q_2(x)+\omega^2 Q_3(x)+O(\omega^3\ln\omega).\\
\end{split}
\end{equation*}
where
\begin{align}\label{Q}
\begin{cases}
\ds Q(x,\omega)=\omega \frac{Im(\varepsilon_c)}{2\pi\gamma_c}\left|u\right|^2,\ms\\
\ds Q_1(x)=\frac{\omega_p^2}{8\pi^3\gamma_c\tau}(\mathcal T_{1}[f])^2,\qquad Q_2(x)=\frac{\omega_p^2}{2\pi^2\gamma_c\tau}\left[\mathcal T_{1}[f]\mathcal N_{D}[f](x)+N(\mathcal T_{1}[f])^2\right],\ms\\
\ds Q_3(x)=\frac{\omega_p^2}{2\pi\gamma_c\tau}\left[ (\mathcal N_{D}[f])^2(x)+\left(N^2+\frac{1}{16}\right)(\mathcal T_{1}[f])^2+2N\mathcal T_{1}[f]\mathcal N_{D}[f](x)\right].\\
\end{cases}
\end{align}

For \eqref{problem2.1}, we assume that the temperature field $v (x, \omega)$ is expressed by
\begin{align}\label{solution V}
v(x,\omega):=
\begin{cases}
\ds V(x)+\mathcal S_D[\psi](x),&\quad x \in \mathbb{R}^{2} \backslash \overline {D}, \ms\\
\ds \mathcal S_D[\phi](x)+\mathcal N_D[Q](x),&\quad x \in D.
\end{cases}
\end{align}
where
\begin{equation}\label{Newton potential}
\begin{split}
\mathcal N_D[Q](x)=&\displaystyle\int_D \Phi_0(x,z)Q(z)dz \\
=&\displaystyle\int_D \Phi_0(x,z)\left[
\omega^2\ln^2\omega Q_1(x)+\omega^2\ln\omega Q_2(x)+\omega^2 Q_3(x)+O(\omega^3\ln\omega)\right]dz,\\
\end{split}
\end{equation}
and the corresponding derivative is
\begin{equation}\label{the derivative of Newton potential}
\begin{split}
\frac{\partial \mathcal N_D[Q](x)}{\partial \nu(x)}=&\displaystyle\int_D \frac{\partial\Phi_0(x,z)Q(z)}{\partial\nu(x)}dz \\
=&\displaystyle\int_D \frac{\partial\Phi_0(x,z)}{\partial\nu(x)}\left[
\omega^2\ln^2\omega Q_1(x)+\omega^2\ln\omega Q_2(x)+\omega^2 Q_3(x)+O(\omega^3\ln\omega)\right] dz.\\
\end{split}
\end{equation}

Using the transmission condition
\begin{equation*}
\begin{cases}
v_+-v_-=0 &\quad \text{on} \ \ \partial D,\\
\frac{\partial v}{\partial \nu}\big{|}_{+}-\gamma_c \frac{\partial v}{\partial \nu}\big{|}_{-}=0 &\quad \text{on} \ \ \partial D.\\
\end{cases}
\end{equation*}
and the jump relation \eqref{eq:jump relation}, we know additionally that
\begin{align*}
\begin{cases}
\ds V(x)+\mathcal S_D[\psi](x)=\mathcal S_D[\phi](x)+\mathcal N_D[Q](x),\ms\\
\ds \frac{\partial H}{\partial \nu}(x)+\left(\mathcal K_D^*+\frac{1}{2}I \right)[\psi](x)=\gamma_c\left(\mathcal K_D^*-\frac{1}{2}I \right)[\psi](x)+\gamma_c\frac{\partial \mathcal N_D[Q]}{\partial \nu}(x),\\
\end{cases}
\end{align*}
that is
\begin{align*}
\begin{cases}
\ds \mathcal S_D[\psi](x)-\mathcal S_D[\phi](x)=\mathcal N_D[Q](x)-V(x),\ms\\
\ds \left(\mathcal K_D^*+\frac{1}{2}I \right)[\psi](x)-\gamma_c\left(\mathcal K_D^*-\frac{1}{2}I \right)=\gamma_c\frac{\partial \mathcal N_D[Q]}{\partial \nu}(x)-\frac{\partial H}{\partial \nu}(x).\\
\end{cases}
\end{align*}
Obviously, we observe that it is possible to convert the above system to the matrix form
\begin{equation*}
\left(
\begin{array}{cc}
  \mathcal S_D & -\mathcal S_D  \\
  \frac{1}{2}I+ \mathcal K_D^*  & -\gamma_c\left(-\frac{1}{2}I+\mathcal K_D^*\right)
\end{array}
\right)
\left(
\begin{array}{c}
   \psi   \\
   \phi
\end{array}
\right)
=
\left(
\begin{array}{c}
  \mathcal N_D[Q]-H  \\
  \gamma_c\frac{\partial \mathcal N_D[Q]}{\partial \nu}-\frac{\partial H}{\partial \nu}
\end{array}
\right).
\end{equation*}
Moreover, we already know the matrix operator
\begin{equation*}
\left(
\begin{array}{cc}
  \mathcal S_D & -\mathcal S_D  \\
  \frac{1}{2}I+ \mathcal K_D^*  & -\gamma_c\left(-\frac{1}{2}I+\mathcal K_D^*\right)
\end{array}
\right)
\end{equation*}
is reversible. Then using Lemma \ref{lem:inverse}, we can obtain
\begin{align*}
\left(
\begin{array}{c}
   \psi   \\
   \phi
\end{array}
\right)
&=
\left(
\begin{array}{cc}
  \mathcal S_D & -\mathcal S_D  \\
  \frac{1}{2}I+ \mathcal K_D^* & -\gamma_c\left(-\frac{1}{2}I+\mathcal K_D^*\right)
\end{array}
\right)^{-1}
\left(
\begin{array}{c}
  \mathcal N_D[Q]-H  \\
  \gamma_c\frac{\partial \mathcal N_D[Q]}{\partial \nu}-\frac{\partial H}{\partial \nu}
\end{array}
\right)\\
&=\frac{1}{\nu_c-1}
\left(
\begin{array}{c}
  \nu_c\left(\lambda I-\mathcal K_D^*\right)^{-1}\left(\frac{1}{2}I-\mathcal K_D^*\right)\mathcal S_D^{-1}\left(\mathcal N_D[Q]-H\right)
  +\left(\lambda I-\mathcal K_D^*\right)^{-1}\left(\gamma_c\frac{\partial \mathcal N_D[Q]}{\partial \nu}-\frac{\partial H}{\partial \nu}\right) \\
  -\left(\lambda I-\mathcal K_D^*\right)^{-1}\left(\frac{1}{2}I+\mathcal K_D^*\right)\mathcal S_D^{-1}\left(\mathcal N_D[Q]-H\right)
  +\left(\lambda I-\mathcal K_D^*\right)^{-1}\left(\gamma_c\frac{\partial \mathcal N_D[Q]}{\partial \nu}-\frac{\partial H}{\partial \nu}\right)
\end{array}
\right).
\end{align*}
Clearly  substituting \eqref{Newton potential} and \eqref{the derivative of Newton potential} into the density $\psi$, we find
\begin{align*}
\psi=&\frac{\gamma_c}{\gamma_c-1}\left(\lambda I-\mathcal K_D^*\right)^{-1}\left(\frac{1}{2}I-\mathcal K_D^*\right)\mathcal S_D^{-1}\left(\mathcal N_D[Q]-H\right)+\frac{1}{\gamma_c-1}\left(\lambda I-\mathcal K_D^*\right)^{-1}\left(\gamma_c\frac{\partial \mathcal N_D[Q]}{\partial \nu}-\frac{\partial H}{\partial \nu}\right)\\
=&\omega^2\ln^2\omega(G_1+F_1)(x)+\omega^2\ln\omega(G_2+F_2)(x)+\omega^2\left(G_3+F_3-\mathcal A_D\right)(x)+O(\omega^3\ln\omega),\\
\end{align*}
where
\begin{equation}\label{def formula2}
\begin{cases}
  \mathcal A_D(x)=\frac{\gamma_c}{\gamma_c-1}(\lambda I-\mathcal K_D^*)^{-1}\left(\frac{1}{2}I- \mathcal K_D^*\right)\mathcal S_D^{-1}V(x)+\frac{1}{\gamma_c-1}(\lambda I-\mathcal K_D^*)^{-1} \frac{\partial V(x)}{\partial \nu},\\
  G_i(x)=\frac{\gamma_c}{\gamma_c-1}\left(\lambda I-\mathcal K_D^*\right)^{-1}\left(\frac{1}{2}I-\mathcal K_D^*\right)\mathcal S_D^{-1}\mathcal N_D[Q_i](x),\\
  F_i(x)=\frac{\gamma_c}{\gamma_c-1}\left(\lambda I-\mathcal K_D^*\right)^{-1}\frac{\partial \mathcal N_D[Q_i](x)}{\partial \nu(x)},\qquad \qquad \qquad \qquad
  i=1,2,3.\\
\end{cases}
\end{equation}
Similarly, by inserting the density $\psi$ into $V(x,\omega)$, we obtain the solution of \eqref{solution V} in $\mathbb{R}^{2} \backslash \overline {D}$:
\begin{equation}\label{V1}
\begin{split}
v(x,\omega)=&V(x)+\mathcal S_D[\psi](x)\\
=&\omega^2\ln^2\omega \mathcal S_D[G_1+F_1](x)+\omega^2\ln\omega \mathcal S_D[G_2+F_2](x)+\omega^2\left\{V(x)+\mathcal S_D\left[G_3+F_3-\mathcal A_D\right](x)\right\}\\
&+O(\omega^3\ln\omega).
\end{split}
\end{equation}
the above identity is necessary to prove the main results.

\section{Proof of Theorem \ref{Th}}\label{sec:proof}
In this section, we devote to the proof of Theorem \ref{Th}.
\begin{proof}[Proof of Theorem \ref{Th}]
\par \textbf{Case 1:} After obtaining the low-frequency asymptotic expansion of the temperature field in Section \ref{sec:Asymptotic analysis-heat}, we can derive integral identities that involve the unknown internal source.

\par Now we proceed with the proof of the main results. Let $v_1(x,\omega)$ and $v_2(x,\omega)$ satisfy the system \eqref{V1} correspond to  $f_1(x) $ and $f_2(x)$. From \eqref{boundary measurement} and \eqref{equation1}, we immediately derive that $\Lambda_{f_1 }(x,\omega)=\Lambda_{f_2  }(x,\omega)$ for $(x,\omega)\in \partial D \times \mathbb{R_+}.$

\par Firstly, we need to compare the coefficient of the term $\omega^2\ln^2\omega$. For any $x\in \partial D$,  we can see that
\begin{equation*}
\mathcal S_D\left[G_1^{(1)} +F_1^{(1)} \right](x)=\mathcal S_D\left[G_1^{(2)} +F_1 ^{(2)} \right](x).
\end{equation*}
Since the single-layer operator $\mathcal S_D$ is injective, we have
\begin{equation*}
\big(G_1^{(1)} +F_1^{(1)} \big)(x)=\big(G_1^{(2)} +F_1^{(2)} \big)(x),
\end{equation*}
and from \eqref{def formula2} it follows that
\begin{equation*}
\begin{split}
&\frac{\gamma_c}{\gamma_c-1}\left(\lambda I-\mathcal K_D^*\right)^{-1}\left(\frac{1}{2}I-\mathcal K_D^*\right)\mathcal S_D^{-1}\mathcal N_D\left[Q_1^{(1)} \right]
+\frac{\gamma_c}{\gamma_c-1}\left(\lambda I-\mathcal K_D^*\right)^{-1}\frac{\partial \mathcal N_D}{\partial \nu}\left[Q_1^{(1)} \right] \\
&=\frac{\gamma_c}{\gamma_c-1}\left(\lambda I-\mathcal K_D^*\right)^{-1}\left(\frac{1}{2}I-\mathcal K_D^*\right)\mathcal S_D^{-1}\mathcal N_D\left[Q_1^{(2)} \right]
+\frac{\gamma_c}{\gamma_c-1}\left(\lambda I-\mathcal K_D^*\right)^{-1}\frac{\partial \mathcal N_D}{\partial \nu}\left[Q_1^{(2)} \right].
\end{split}
\end{equation*}
We then combine \eqref{def formula1} and \eqref{Q} to discover
\begin{equation}\label{f1}
\begin{split}
&\left(\frac{1}{2}I-\mathcal K_D^*\right)\mathcal S_D^{-1} \mathcal N_D\left[\left(\mathcal T_{1}[f]^{(1)} \right)^2\right]+\frac{\partial \mathcal N_D}{\partial \nu}\left[\left(\mathcal T_{1}[f]^{(1)} \right)^2\right]\\
&=\left(\frac{1}{2}I-\mathcal K_D^*\right)\mathcal S_D^{-1}\mathcal N_D\left[\left(\mathcal T_{1}[f]^{(2)} \right)^2\right]+\frac{\partial \mathcal N_D}{\partial \nu}\left[\left(\mathcal T_{1}[f]^{(2)} \right)^2\right],
\end{split}
\end{equation}
where $$ \mathcal T_{1}[f]^{(1)} =\displaystyle\int_Df_1(y)dy, \quad \mathcal T_{1}[f]^{(2)} =\displaystyle\int_Df_2(y)dy. $$
Since $\mathcal T_{1}[f]$ is constant, and consequently \eqref{f1} implies
\begin{equation*}
\left[\left(\frac{1}{2}I-\mathcal K_D^*\right)\mathcal S_D^{-1}\mathcal N_D[1](x)+\frac{\partial \mathcal N_D}{\partial \nu}[1](x)\right]\left[\left(\mathcal T_{1}[f]^{(1)} \right)^2-\left(\mathcal T_{1}[f]^{(2)} \right)^2\right]=0.\\
\end{equation*}
Our task next is to analyze the above equation. Indeed, we first note that
\begin{align*}
\begin{cases}
\ds \Delta (\mathcal N_D[1])(x)=1,&\qquad x\in D,\ms\\
\ds \Delta (\mathcal N_D[1])(x)=0,&\qquad x\in \mathbb{R}^2\backslash\overline{D}, \\
\end{cases}
\end{align*}
According to Lemma \ref{lem:inverse}, it follows that
\begin{equation}\label{B}
\begin{split}
&\left[\left(\frac{1}{2}I-\mathcal K_D^*\right)\mathcal S_D^{-1}\mathcal N_D[1](x)+\frac{\partial \mathcal N_D}{\partial \nu}[1](x)\right]\left[\left(\mathcal T_{1}[f]^{(1)} \right)^2-\left(\mathcal T_{1}[f]^{(2)} \right)^2\right]\\
&=T_f^{(1)} (x)\left[\left(\mathcal T_{1}[f]^{(1)} \right)^2-\left(\mathcal T_{1}[f]^{(2)} \right)^2\right],
\end{split}
\end{equation}
where
\begin{equation*}
\begin{split}
T_f^{(1)} (x)&=\left(\int_D 1dx-\int_{\partial D}\mathcal S_D^{-1}\big[\mathcal N_D[1]\big](y)d\sigma(y) \right)\varphi_0+\mathcal S_D^{-1}\big[\mathcal N_D[1]\big]\\
&=\left(|D|-\int_{\partial D}\mathcal S_D^{-1}\mathcal N_D[1](y)d\sigma(y)\right)\varphi_0+\mathcal S_D^{-1}\mathcal N_D[1].
\end{split}
\end{equation*}
Since $ \mathcal N_D[1]\neq 0,$ we discover that
\begin{equation*}
\left(|D|-\int_{\partial D}\mathcal S_D^{-1}\mathcal N_D[1](y)d\sigma(y)\right)\varphi_0\in {c_0\varphi_0},
\end{equation*}
where $c_0$ is an constant. In addition
\begin{equation*}
\mathcal S_D^{-1}\mathcal N_D[1]\in H_0^*(\partial D)=\{c_0\varphi_0 \}^{\perp},
\end{equation*}
we thus deduce that
\begin{equation*}
T_f^{(1)} (x)\in {c_0\varphi_0}\oplus H_0^*(\partial D)\neq 0.
\end{equation*}
Then by combining with \eqref{B}, we consequently have the phaseless integral identity of internal source $f$,
\begin{equation}\label{r1}
\left|\int_D f_1(y)dy\right|=\left|\int_D f_2(y)dy\right|.
\end{equation}

\par Secondly, we also need to compare the coefficient of the term $\ln\omega$. By a similar argument, we can see that for any $x\in \partial D,$
\begin{equation*}
\begin{split}
&\left(\frac{1}{2}I-\mathcal K_D^*\right)\mathcal S_D^{-1}\mathcal N_D\left[\mathcal N_D[f_1]\right]\mathcal T_{1}[f]^{(1)} +\frac{\partial \mathcal N_D\left[\mathcal N_Df_1\right]}{\partial \nu(x)}\mathcal T_{1}[f]^{(1)}  \\
&=\left(\frac{1}{2}I-\mathcal K_D^*\right)\mathcal S_D^{-1}\mathcal N_D\left[\mathcal N_D[f_2]\right]\mathcal T_{1}[f]^{(2)} +\frac{\partial \mathcal N_D\left[\mathcal N_Df_2\right]}{\partial \nu(x)}\mathcal T_{1}[f]^{(2)},\\
\end{split}
\end{equation*}
that is
\begin{equation*}
\left[\left(\frac{1}{2}I-K_D^*\right)\mathcal S_D^{-1}\mathcal N_D+\frac{\partial}{\partial \nu}\mathcal N_D\right]\big[\mathcal T_{1}[f]^{(1)} \mathcal N_D[f_1]-\mathcal T_{1}[f]^{(2)} \mathcal N_D[f_2]\big](x)=0.\\
\end{equation*}
Note that the volume potential
\begin{align*}
\mathcal N_D[\mathcal N_D[f]](x)&=\int_D\Phi_0(x,z)\int_D\Phi_0(z,y)f(y)dydz\\
&=\int_D\int_D\Phi_0(x,z)\Phi_0(z,y)f(y)dydz,
\end{align*}
satisfies
\begin{align*}
\begin{cases}
\ds \Delta (\mathcal N_D[\mathcal N_D[f]])(x)=\mathcal N_D[f](x),&\quad x\in D,\ms\\
\ds \Delta (\mathcal N_D[\mathcal N_D[f]])(x)=0,&\quad x\in \mathbb{R}^2\backslash\overline{D}.\\
\end{cases}
\end{align*}
According to Lemma \ref{lem:inverse}, it follows that
\begin{align*}
\left(\left(\frac{1}{2}I-\mathcal K_D^*\right)\mathcal S_D^{-1}\mathcal N_D[\mathcal N_D[f]](x)+\frac{\partial}{\partial \nu}\mathcal N_D[\mathcal N_D[f]](x)\right)\mathcal T_{1}[f]:=T_f^{(2)} (x) \mathcal T_{1}[f],
\end{align*}
where
\begin{equation*}
\begin{split}
T_f^{(2)} (x)&=\left(\int_D \mathcal N_D[f](x)dx-\int_{\partial D}\mathcal S_D^{-1}\big[\mathcal N_D[\mathcal N_D[f]]\big](y)d\sigma(y) \right)\varphi_0(x)+\mathcal S_D^{-1}\big[\mathcal N_D^2[f]\big](x).
\end{split}
\end{equation*}
Then, arguing as above, we likewise deduce
\begin{align*}
T_f^{(2)} (x) \in {\tilde{c}_0\varphi_0}\oplus H_0^*(\partial D),
\end{align*}
where $\tilde{c}_0$ is an constant. We can further compute that
\begin{equation*}
 \mathcal T_{1}[f]^{(1)} \mathcal S_D^{-1}\mathcal N_D^2[f_1](x)=\mathcal T_{1}[f]^{(2)} \mathcal S_D^{-1}\mathcal N_D^2[f_2](x) ,\qquad x\in \partial D.
\end{equation*}
In view of the identity \eqref{r1}, we obtain that
\begin{equation*}
\mathcal N_D^2[f_1](x)=\pm \mathcal N_D^2[f_2](x),\qquad x\in \partial D.
\end{equation*}

Owing to the properties of volume potential, we obtain
\begin{align*}
\begin{cases}
\ds  \Delta^2 \big(\mathcal N_D^2[f_i](x)\big)=f_i(x),&\qquad x\in D,\ms\\
\ds  \Delta \big(\mathcal N_D^2[f_i](x)\big)=0, &\qquad x\in\mathbb{R}^{2} \backslash \overline {D}.
\end{cases}
\end{align*}
Hence we can apply the uniqueness of the exterior Dirichlet problem to the harmonic equation, and we derive
\begin{equation*}
\mathcal N_D^2[f_1](x)=\pm \mathcal N_D^2[f_2](x),\qquad x\in \mathbb{R}^{2} \backslash \overline {D}.
\end{equation*}
%
Setting $\tilde{w}=\mathcal N_D^2[f_1\pm f_2](x)$, $\tilde{f}=f_1\pm f_2$, and notice that $\tilde{w}=0$ in $\mathbb{R}^{2} \backslash \overline {D}$, we obtain the following inhomogeneous biharmonic system with Navier boundary conditions \cite{Gazz2010}:
\begin{align*}
\begin{cases}
\ds    \Delta^2 \tilde{w}=\tilde{f}  &\quad \text{in} \ \  D,\ms\\
\ds    \tilde{w}=0  &\quad \text{on} \ \ \partial D,\ms\\
\ds    \Delta\tilde{w}=0  &\quad \text{on} \ \ \partial D.
\end{cases}
\end{align*}
By Green's formula, for arbitrary biharmonic function $\tilde{v}\in H^4(D)$, it follows that
\begin{equation}\label{green}
\int_D\tilde{f}\tilde{v}dx=\int_D\tilde{v}\Delta^2\tilde{w}dx-\int_D\tilde{w}\Delta^2\tilde{v}dx
=\int_{\partial D}\left(\tilde{v}\frac{\partial\Delta\tilde{w}}{\partial \nu}+\Delta\tilde{v}\frac{\partial\tilde{w}}{\partial \nu}-\tilde{w}\frac{\partial\Delta\tilde{v}}{\partial \nu}-\Delta\tilde{w}\frac{\partial\tilde{v}}{\partial \nu}\right)d\sigma(x)=0.
\end{equation}
Since $f_1$ and $f_2$ differ by a harmonic part or biharmonic part, by taking $\tilde{v}=\tilde{f}$ in (\ref{green}), we get
\begin{equation*}
\int_D|\tilde{f}|^2dx=0,
\end{equation*}
and the phaseless uniqueness holds.

\textbf{Case 2:} Owing to the rotation invariance of the biharmonic operator $\Delta^2$, the vector $\eta$ can be suitably rotated to any coordinate axis. Without loss of generality, we assume that $\eta$ is rotated to align with the $x_2$-axis. Consequently, it implies that $\partial_{x_2}\tilde{f}=0$.

Setting $\tilde{v}=e^{i\zeta\cdot x}$, $x\in\mathbb{R}^{2}$, we have
\begin{equation*}
\zeta=\zeta_1+i\zeta_2\in\mathbb{C}^{2},\ \ \zeta_1=(\zeta_{11},0)\in\mathbb{R}^{2},\ \ \zeta_1=(0,\zeta_{22})\in\mathbb{R}^{2},\ \ \text{and}\ \ \zeta_{11} ^2-\zeta_{22}^2=0.
\end{equation*}
It is easy to check that $e^{i\zeta\cdot x}$ is a harmonic function and also a biharmonic function in $\mathbb{R}^{2}$. Plugging $\tilde{v}=e^{i\zeta\cdot x}$ into (\ref{green}), we find
\begin{align*}
\int_D\tilde{f}\tilde{v}dx=\int_{\mathbb{R}^{2}}\tilde{f}(x_1)e^{i\zeta\cdot x}dx=\int_{\mathbb{R}^{1}}\tilde{f}(x_1)e^{i\zeta_{11}\cdot x_1}dx_1\cdot
\int_{\{x_2:(x_1,x_2)\in B_R\}}e^{-i\zeta_{22}\cdot x_2}dx_2=0.
\end{align*}
Then $\tilde{f}=0$, and
\begin{equation*}
 \big|f_1(x)\big|=\big| f_2(x)\big|,\qquad x\in D. \\
\end{equation*}
The proof is completed.
\end{proof}

\section*{Acknowledgements}
\addcontentsline{toc}{section}{Acknowledgments}
The research of H. Liu was supported by NSFC/RGC Joint Research Scheme, N CityU101/21, ANR/RGC
Joint Research Scheme, A-CityU203/19, and the Hong Kong RGC General Research Funds (projects 11311122, 11304224 and 11300821).
The research of G. Zheng was supported by the NSF of China (12271151), NSF of
Hunan (2020JJ4166) and NSF Innovation Platform Open Fund project of Hunan Province (20K030).

\end{document}